\documentclass{article}
\usepackage{amssymb}
\usepackage{amsmath}
\usepackage{cite}

\newtheorem{theorem}{Theorem}

\newtheorem{corollary}[theorem]{Corollary}

\input{tcilatex}
\begin{document}

\title{Lie symmetries for systems of evolution equations}
\author{Andronikos Paliathanasis$^{1,2}$\thanks{%
Email: anpaliat@phys.uoa.gr} \\
{\ }$^{1}${\textit{Instituto de Ciencias F\'{\i}sicas y Matem\'{a}ticas, }}\\
{\textit{Universidad Austral de Chile, Valdivia, Chile}}\\
$^{2}${\textit{Institute of Systems Science, Durban University of Technology}%
}\\
{\textit{Durban 4000, Republic of South Africa}} \and Michael Tsamparlis$%
^{3} $\thanks{%
Email: mtsampa@phys.uoa.gr} \\
{\ }$^{3}${\textit{Faculty of Physics, Department of
Astronomy-Astrophysics-Mechanics,}}\\
{\ \textit{University of Athens, Panepistemiopolis, Athens 157 83, Greece}}}
\maketitle

\begin{abstract}
The Lie symmetries for a class of systems of evolution equations are
studied. The evolution equations are defined in a bimetric space with two
Riemannian metrics corresponding to the space of the independent and
dependent variables of the differential equations. The exact relation of the
Lie symmetries with the collineations of the bimetric space is determined.

\end{abstract}

\section{Introduction}

Lie symmetries is a powerful method for the determination of solutions in
the theory of differential equations. A Lie symmetry is important as it
provides invariants which can be used to write a new differential equation
with less degree of freedom. Furthermore, a solution of such an equation is
also a solution of the original differential equation \cite{Stephani,Bluman}%
. The reduction process is the main application of Lie symmetries, however,
it is not a univocal approach. Symmetries can be used for the determination
of conservation currents \cite{noe}, for the classification of differential
equations \cite{ovsiannikov,cla1,cla2,cla3,cla4,cla5,clas6} and for the
reconnaissance of some well-known systems \cite%
{leach1,leach2,leach3,leach4,erm}.

In the recent literature, it has been shown that there is a close relation
between the Lie symmetries of a second order differential equation and the
geometry of the space where motion occurs. For example, the conservation of
energy and angular momentum in Newtonian Physics is a result of the Lie
point symmetries, generated by the Killing vectors of translations and
rotations respectively. The general result for a holonomic autonomous
dynamical system moving in a Riemannian space is that the Lie point
symmetries of the equations of motion are generated by the special
projective algebra of the kinematic metric, that is, the metric defined by
the kinetic energy. \cite{katz1,katz2,katz3,ami1,ami2,p1,p2}. Similar
results have been derived for second-order partial differential equations 
\cite{p3,p4,pp05}.

The importance of the relation of Lie symmetries with the symmetry groups
lies in that one can determine those symmetries, hence, reduce the problem
of finding solutions for second order time independent differential
equations by applying direct theorems of differential geometry \cite%
{p6,p7,p8}. Moreover, it has been shown \cite{chris} that the application of
the geometric approach in constrained dynamical systems can be used to
extract information during the quantization process\footnote{%
For some interesting results on the quantization from symmetries see \cite%
{nuz1,nuz2,nuz4} and references therein.} \cite{chris2,chris3,chris4,anp1}.

The present work is a natural extension of the approach initiated in \cite%
{p3,p4}, that is, we consider the problem of determining the Lie (point)
symmetries for a class of systems of quasilinear evolution equations which
define a bimetric manifold with two independent connections. The system of
differential equations is assumed to be of the form%
\begin{equation}
Q^{A}\left( t,x^{k},u^{C},u_{,i}^{C},u_{,ij}^{C},u_{t}^{C}\right) =0,
\label{ht.01}
\end{equation}%
in which%
\begin{equation}
Q^{A}\equiv g^{ij}u_{,ij}^{A}+g^{ij}\tilde{\Gamma}%
_{BC}^{A}u_{,i}^{B}u_{,j}^{C}-\Gamma ^{i}u_{,i}^{A}+F^{A}\left(
t,x^{k},u^{C}\right) -u_{,t}^{A},  \label{ht.02}
\end{equation}%
where $u^{A}$ denotes the $m-$dependent variables, $\dim \left( u^{A}\right)
=m$, $x^{k}$ are the $n-$independent variables, $\dim \left( x^{k}\right) =n$%
, while $\Gamma ^{i}=g^{jk}\Gamma _{jk}^{i}\left( x^{r}\right) $ and $\tilde{%
\Gamma}_{BC}^{A}=\tilde{\Gamma}_{BC}^{A}\left( u^{D}\right) $ are the two
affine connections on the bimetric theory $\left\{ g_{ij},~H_{AB}\right\} $.
Metric $g_{ij}$ is defined in the space of the independent variables while $%
H_{AB}\left( u^{C}\right) $ defines the geometry of the space of dependent
variables $u^{A}$. Moreover, we assume that the two spaces are minimally
coupled, that is, $\left\{ g_{ij},~H_{AB}\right\} =\left\{ g_{ij}\left(
x^{k}\right) ,H_{AB}\left( u^{C}\right) \right\} .$ Finally, the function $%
F^{A}\left( t,x^{k},u^{C}\right) $ counts for the interaction term of the
metrics $g_{ij},~H_{AB}$ and $t$ denotes the \textquotedblleft
time\textquotedblright\ in which the system evolves.

Static solutions of the system (\ref{ht.01}) have been studied in \cite{p4}.
It was found that the Lie point symmetries of (\ref{ht.01}) are constructed
by the conformal Killing vectors and the affine collineations of the metrics 
$g_{ij}$ and $H_{AB}$ respectively. The present letter is organized as
follows.

In Section \ref{sec2} we present the main results. We solve the symmetry
conditions for the system (\ref{ht.01}) and show that the Lie point
symmetries are generated by the elements of the homothetic algebra and the
affine algebra of the two metrics which are defined by the differential
equations.\ Three corollaries are presented for some cases of special
interest. Finally, our discussion is given in Section \ref{conclusions}.

\section{Lie point symmetries}

\label{sec2}

The lhs of the differential equation (\ref{ht.01}) can be seen as a vector
field $Q^{A}$ in the space of variables $A=A\left(
t,x^{k},u^{C},u_{,i}^{C},u_{,ij}^{C},u_{t}^{C}\right) .~$ The condition is
that the function $Q^{A}$ is invariant under the action of the one parameter
point transformation 
\begin{align}
\bar{t}& =t+\varepsilon \xi ^{t}\left( t,x^{k},u^{B}\right)  \label{ht.03} \\
\bar{x}^{i}& =x^{i}+\varepsilon \xi ^{i}\left( t,x^{k},u^{B}\right) ,
\label{ht.04} \\
\bar{u}^{A}& =\bar{u}^{A}+\varepsilon \eta ^{A}\left( t,x^{k},u^{B}\right) ~,
\label{ht.05}
\end{align}%
is 
\begin{equation}
\mathbf{X}^{[2]}Q^{A}=0  \label{ht.06}
\end{equation}%
where $\mathbf{X}^{[2]}$ is the second extension for the generator $X=\xi
^{t}\partial _{t}+\xi ^{i}\partial _{i}+\eta ^{A}\partial _{A}$ in the space
of both independent and dependent variables. In particular, the second
extension is given by the formula \cite{Stephani,Bluman}%
\begin{equation}
\mathbf{X}^{[2]}=\mathbf{X}+\eta _{\alpha }^{A}\partial _{u_{\alpha
}^{A}}+\eta _{\alpha \beta }^{A}\partial _{u_{\alpha \beta }^{A}}
\label{ht.07}
\end{equation}%
where\footnote{%
The Greek indices~$a,\beta ,\gamma $ count on all the independent parameters~%
$x^{i}$ and $t.$}~$\eta _{\alpha }^{A}$ and~$\eta _{\alpha \beta }^{A}$ are
defined as follows 
\begin{equation}
\eta _{\alpha }^{A}=\eta _{,\alpha }^{A}+u_{,\alpha }^{B}\eta _{,B}^{A}-\xi
_{,\alpha }^{\beta }u_{,\beta }^{A}-u_{,\beta }^{A}u_{,\alpha }^{B}\xi
_{,B}^{\beta }~,  \label{ht.08}
\end{equation}%
\begin{align}
\eta _{\alpha \beta }^{A}& =\eta _{,\alpha \beta }^{A}+2\eta _{,B(\alpha
}^{A}u_{,\beta )}^{B}-\xi _{,\alpha \beta }^{\gamma }u_{,\gamma }^{A}+\eta
_{,BC}^{A}u_{,\alpha }^{B}u_{,\beta }^{C}-2\xi _{,(\alpha |B|}^{\gamma
}u_{\beta )}^{B}u_{,\gamma }^{A}  \notag \\
& -\xi _{,BC}^{\gamma }u_{,\alpha }^{B}u_{,\beta }^{C}u_{,\gamma }^{A}+\eta
_{,B}^{A}u_{,\alpha \beta }^{B}-2\xi _{,(\beta }^{\gamma }u_{,\alpha )\gamma
}^{A}-\xi _{,B}^{\gamma }\left( u_{,\gamma }^{A}u_{,\alpha \beta
}^{B}+2u_{(,\beta }^{B}u_{,\alpha )\gamma }^{A}\right) .  \label{ht.09}
\end{align}

Considering condition, (\ref{ht.06}), as a polynomial in the variables ~$%
\left( u^{A}\right) ^{0}$,~$\left( u_{,\alpha }^{A}\right) ,~\left(
u_{,\alpha }^{A}u_{,\beta }^{B}\right) $,$\left( u_{,\alpha \beta
}^{A}\right) $ and$~\left( u_{,\alpha \beta }^{A}u_{,\gamma }^{B}\right) $
we find the following symmetry conditions 
\begin{equation}
\xi _{,i}^{t}=0~,~\xi _{,A}^{t}=0~,~\xi _{,A}^{i}=0  \label{ht.10}
\end{equation}%
\begin{equation}
F_{,k}^{A}\xi ^{k}+F_{,t}^{A}\xi ^{t}+F_{,B}^{A}\eta ^{B}+\left( g^{ij}\eta
_{;ij}^{A}-\eta _{,t}^{A}\right) =\lambda F^{A}  \label{ht.11}
\end{equation}%
\begin{equation}
\left( -g^{kj}\xi _{,kj}^{i}-\Gamma _{~,k}^{i}\xi ^{k}+\xi _{,j}^{i}\Gamma
^{j}+\lambda \Gamma ^{i}+\xi _{,t}^{i}\right) \delta _{B}^{A}-\Gamma
^{i}\eta _{,B}^{A}+2g^{ij}\left( \eta _{|B}^{A}\right) _{,j}=0  \label{ht.12}
\end{equation}%
\begin{equation}
\left( L_{\xi }g^{ij}-\lambda g^{ij}\right) \tilde{\Gamma}_{BC}^{A}+\tilde{%
\Gamma}_{BC,D}^{A}\eta ^{D}+g^{ij}\eta _{,BC}^{A}+2g^{ij}\tilde{\Gamma}%
_{KC}^{A}\eta _{,C}^{K}=0  \label{ht.13}
\end{equation}%
\begin{equation}
\delta _{B}^{A}L_{\xi }g^{ij}=\left( \lambda -\eta _{,B}^{A}\right)
g^{ij}\delta _{B}^{A}  \label{ht.14}
\end{equation}%
\begin{equation}
\eta _{,B}^{A}=\left( \xi _{,t}^{t}-\lambda \right) \delta _{B}^{A}
\label{ht.14a}
\end{equation}%
where \textquotedblleft $|$\textquotedblright\ and \textquotedblleft $;$%
\textquotedblright\ denote covariant derivative with respect to the
connection $\tilde{\Gamma}_{BC}^{A}$ or $\Gamma _{jk}^{i}$ respectively,
while $\mathcal{L}_{X}$ denotes Lie derivative. The symmetry conditions are
similar to those derived in \cite{p4}, therefore, we follow the same
procedure in their solution.

Condition (\ref{ht.14}) becomes 
\begin{equation}
L_{\xi }g_{ij}=\left( \lambda -\eta _{,B}^{A}\right) g_{ij}  \label{ht.15a}
\end{equation}%
which means that $\xi ^{i}$ is a conformal Killing vector (CKV) field for
the metric $g_{ij}$ with conformal factor $\psi =\frac{1}{2}\left( \lambda
-\eta _{,B}^{A}\right) $. However, from (\ref{ht.14a}) we have that $\psi
=\psi \left( t\right) $ only such that $\xi ^{t}=2\psi \int T\left( t\right)
dt.~$Hence, $\xi ^{i}=T\left( t\right) \zeta ^{i},~$where $\zeta ^{i}$\ is a
homothetic vector (HV) of $g_{ij}.$

From (\ref{ht.13}) we have that%
\begin{equation}
\mathcal{L}_{\eta }\tilde{\Gamma}_{BC}^{A}=0  \label{ht.15}
\end{equation}%
which implies that $\eta ^{A}$ is an affine collineation (AC) of the metric $%
H_{AB}$. \ These two conditions, namely (\ref{ht.15a}) and (\ref{ht.15})
provide the geometric answer to the relation between Lie point symmetries of
(\ref{ht.01}) and collineations of the underlying space.

The remaining symmetry conditions produce constraints on these collineations
which should be satisfied in order the latter to be Lie point symmetries for
the system (\ref{ht.01}).

We rewrite condition (\ref{ht.12}) as follows%
\begin{equation}
\left( \eta _{|B}^{A}\right) _{,j}=-T_{,t}\zeta _{i}\delta _{B}^{A}
\label{ht.16a}
\end{equation}%
from which follows that 
\begin{equation}
\eta _{|B}^{A}=-T_{,t}\bar{\zeta}\delta _{B}^{A}+Z\left( t,u^{C}\right)
\label{ht.17b}
\end{equation}%
where $\bar{\zeta}_{,i}$ denotes a gradient HV of $g_{ij}$, otherwise $%
T\left( t\right) =const.$ such that $T_{,t}=0$. From (\ref{ht.17b}) and (\ref%
{ht.15}) it follows that%
\begin{equation}
\eta ^{A}=-T_{,t}\bar{\zeta}Y^{A}\left( u^{C}\right) +Z\left( t,u^{C}\right)
+\sigma ^{A}\left( t,x^{k}\right)  \label{ht.18}
\end{equation}%
where $Z^{A}$ is an AC of $H_{AB}$ and $Y^{A}$ is a proper gradient HV of $%
H_{AB}.$

Finally, condition (\ref{ht.11}) becomes%
\begin{equation}
\mathcal{L}_{\xi ^{i}}F^{A}+2\psi \int T\left( t\right) dtF_{,t}^{A}+2\psi
T\left( t\right) F^{A}-T_{,t}\bar{\zeta}\mathcal{L}_{Y}F^{A}+\mathcal{L}%
_{Z}F^{A}+\left( g^{ij}\eta _{;ij}-\eta _{,t}\right) =0  \label{ht.19}
\end{equation}

The last condition is the one which constrains the collineations of the
bimetric manifold with the function $F^{A}$ in order for the collineations
to generate a Lie point symmetry for the system (\ref{ht.01}). The following
theorem is the main result of this work and a natural extension of \cite{p4}
for systems of evolution equations.

\begin{theorem}
\label{Lsym}The Lie point symmetries for the system of evolution equations (%
\ref{ht.01}) are generated by the collineations of the bimetric $\left\{
g_{ij},~H_{AB}\right\} .~$In particular, the generic form for the Lie
symmetry vector is given by the formula%
\begin{equation}
X=2\psi \int T\left( t\right) dt\partial _{t}+T\left( t\right) \zeta
^{i}\left( x^{k}\right) \partial _{i}+\left( Z^{A}\left( t,u^{A}\right)
-T_{,t}\bar{\zeta}Y^{A}\left( u^{C}\right) +\sigma ^{A}\left( t,x^{k}\right)
\right) \partial _{A}  \label{ht.20}
\end{equation}%
where $\zeta ^{i}$ is an element of the Homothetic algebra of $g_{ij}~$with
homothetic factor $\psi ,~\bar{\zeta}$ is a gradient homothetic vector field
of $g_{ij}$, $Z^{A}$ and $Y^{A}$ are elements of the affine and homothetic
algebra of the metric $H_{AB}$ respectively. The collineations and the
functions $T\left( t\right) ,~\sigma \left( t,x^{k}\right) $ are constrained
by condition (\ref{ht.19}).
\end{theorem}

In the case where the bimetric does not admit any gradient HV, it follows:

\begin{corollary}
\label{Lsym1}The generic Lie point symmetry for the system of evolution
equations (\ref{ht.01}) in which $g_{ij}$ and $H_{AB}$ do not admit a
gradient HV is%
\begin{equation}
X=2\psi t\partial _{t}+\zeta ^{i}\left( x^{k}\right) \partial _{i}+\left(
Z^{A}\left( t,u^{A}\right) +\sigma ^{A}\left( t,x^{k}\right) \right)
\partial _{A},  \label{ht.21}
\end{equation}%
where $\zeta ^{i}$,~$Z^{A}$ and $\sigma \left( t,x^{k}\right) $ are
constrained with the condition%
\begin{equation}
\mathcal{L}_{\xi ^{i}}F^{A}+2\psi tF_{,t}^{A}+2\psi F^{A}+\mathcal{L}%
_{Z}F^{A}+\left( g^{ij}\sigma _{;ij}^{A}-\sigma _{,t}^{A}\right)
-Z_{,t}^{A}=0.  \label{ht.22}
\end{equation}
\end{corollary}

In case where there is no interaction term $F^{A}$ in the system (\ref{ht.01}%
), that is $F^{A}\left( t,x^{k},u^{C}\right) =0,$ condition (\ref{ht.19})
becomes%
\begin{equation}
\left( g^{ij}\sigma _{;ii}^{A}-\sigma _{,t}^{A}\right) +\left( \psi
T_{,t}-T_{,tt}\bar{\zeta}\right) Y^{A}-Z_{,t}^{A}=0  \label{ht.23}
\end{equation}%
that is 
\begin{equation}
\psi T_{,t}Y^{A}-Z_{,t}^{A}=0~,~T_{,tt}\bar{\zeta}Y^{A}=0,  \label{ht.24}
\end{equation}%
and%
\begin{equation}
g^{ij}\sigma _{;ii}^{A}-\sigma _{,t}^{A}=0.  \label{ht.25}
\end{equation}

From (\ref{ht.24}) it follows that $T_{,tt}=0$, that is $T\left( t\right)
=T_{0}+T_{1}t$, while from (\ref{ht.24}) we have that~$Z^{A}=\psi
T_{1}tY^{A} $ or $Z^{A}\neq Y^{A}$ and $T_{1}$. Condition (\ref{ht.25})
provides an infinite number of trivial symmetries, which is a well-known
result for the heat equation. Thus, from Theorem \ref{Lsym} follows:

\begin{corollary}
\label{Lsym2}The generic Lie point symmetry for the system of free evolution
equations 
\begin{equation}
g^{ij}u_{,ij}^{A}+g^{ij}\tilde{\Gamma}_{BC}^{A}u_{,i}^{B}u_{,j}^{C}-\Gamma
^{i}u_{,i}^{A}+F^{A}\left( t,x^{k},u^{C}\right) -u_{,t}^{A}  \label{ht.26}
\end{equation}%
is generated by the HV of the metric $g_{ij}$ and the ACs of the connection $%
\tilde{\Gamma}_{BC}^{A}$ and it is given by the following formula%
\begin{eqnarray}
X &=&\left( \bar{\psi}\left( T_{1}t^{2}+T_{0}t\right) +\alpha _{1}\psi
t+\alpha _{0}\right) \partial _{t}+\left( \bar{\psi}\left(
T_{1}t+T_{0}\right) \bar{\zeta}^{,i}\left( x^{k}\right) +\alpha _{1}\zeta
^{i}\right) \partial _{i}+  \notag \\
&&~~~~~+\left( Z^{A}\left( t,u^{A}\right) -\bar{\psi}T_{1}\left( t+1\right) 
\bar{\zeta}Y^{A}\left( u^{C}\right) +\sigma ^{A}\left( t,x^{k}\right)
\right) \partial _{A}  \label{ht.27}
\end{eqnarray}%
in which $\zeta ^{i},\bar{\zeta}^{,i}~$are elements of the Homothetic
algebra of $g_{ij}$ with Homothetic factors $\psi $ and~$\bar{\psi}$, where $%
\bar{\zeta}^{,i}$ is a gradient vector field. $Z^{A}$ is an affine
collineation of $\tilde{\Gamma}_{BC}^{A}$ while $Y^{A}$ is a HV for the
metric $H_{AB}$. Finally, $\sigma ^{A}\left( t,x^{k}\right) $ solves the
system of independent linear evolution equations (\ref{ht.25}).
\end{corollary}

In the special scenario in which $\left\{ g_{ij},~H_{AB}\right\} $ are
metrics of maximally symmetric spaces with nonzero curvature from
Corollaries \ref{Lsym1} and \ref{Lsym2}, we have:

\begin{corollary}
\label{Lsym3}The generic Lie point symmetry for the system of free evolution
equations of the form (\ref{ht.26}) where $\left\{ g_{ij},~H_{AB}\right\} $
are metrics of maximally symmetric spaces with nonzero curvature is given as
below%
\begin{equation}
X=\alpha _{0}\partial _{t}+\zeta ^{i}\partial _{i}+\left( Z^{A}\left(
t,u^{A}\right) +\sigma ^{A}\left( t,x^{k}\right) \right) \partial _{A}
\end{equation}%
where $\zeta ^{i}$ and $Z^{A}$ are the $\frac{n\left( n+1\right) }{2}$ and $%
\frac{m\left( m+1\right) }{2}$ Killing vectors of the spaces $\left\{
g_{ij},~H_{AB}\right\} $. Finally, $\sigma ^{A}\left( t,x^{k}\right) $ solve
the system of independent linear evolution equations (\ref{ht.25}).
\end{corollary}

For the proof of Corollary \ref{Lsym3}, we remind to the reader that
maximally symmetric spaces of nonzero curvature do not admit a HV and ACs. \ 

At this point, we omit the demonstration of the application of the derived
results. The same approach which has been presented in \cite{p4} can be
applied.

\section{Conclusion}

\label{conclusions}

This work focuses on the relation between geometry and algebraic properties
of differential equations. The system of evolution equations (\ref{ht.01})
is defined by a bimetric theory where the two Riemannian spaces define the
spaces of the independent and dependent variables.

The main result of this work is Theorem \ref{Lsym}, in which we construct
the generic form of the Lie point symmetries for the system (\ref{ht.01})
and we transfer the problem of the determination of point Lie symmetries to
a geometric problem. In particular, the geometric construction of the
bimetric theory defines the existence, the number and the form of the Lie
point symmetries.

Because the point transformation (\ref{ht.03})-(\ref{ht.05}) has components
in the spaces of the independent and the dependent variables, it follows
that it is constructed by the HV and the ACs for the two metrics which
define the bimetric theory. \ There are various differences with the
quasilinear system we studied in \cite{p4} the main one being that in \cite%
{p4} the CKVs of the metric $g_{ab}$ generate symmetries of the system. That
means that in the case of static solutions of the system (\ref{ht.01}) the
new \textquotedblleft Type-II\textquotedblright\ hidden symmetries \cite%
{pp1,pp2} are generated by the proper CKVs of the metric $g_{ij}$, a result
which generalizes the corresponding result of the one-dimensional case \cite%
{pp3}.

With the present work we complete our analysis on the relation between
geometry and Lie point symmetries of differential equations of second-order
of physical interest. There are various ways through which one could extend
that geometric approach by studying other families of differential equations
or studying higher-order symmetries. A review on the geometric approach of
the symmetries of differential equations where some extensions will be
presented is in progress.

\bigskip%

{\large {\textbf{Acknowledgements}}} \newline
The research of AP was supported by FONDECYT postdoctoral grant no. 3160121.
AP thanks the University of Athens for the hospitality provided while part
of this work took place.


\begin{thebibliography}{99}
\bibitem{Stephani} H. Stephani, Differential Equations: Their Solutions
Using Symmetry, Cambridge University Press, New York, (1989)

\bibitem{Bluman} G.W. Bluman and S. Kumei, Symmetries of Differential
Equations, Springer-Verlag, New York, (1989)

\bibitem{noe} E. Noether, Nachr. v.d. Ges. d. Wiss. zu Gottingen \textbf{235}%
, (1918)

\bibitem{ovsiannikov} L. V. Ovsiannikov, \ Group analysis of differential
equations, Academic Press, New York, (1982)

\bibitem{cla1} N. Kallinikos and E. Meletidou, J. Phys. A: Math. Theor. 
\textbf{46,} 305202 (2013)

\bibitem{cla2} N. Kallinikos, Group classification of charged particle
motion in stationary electromagnetic fields, arXiv:1707.05684

\bibitem{cla3} K. Andriopoulos, P.G.L. Leach and G.P. Flessas, J. Math.
Anal. Appl. \textbf{262,} 256 (2001)

\bibitem{cla4} A.H. Kara and F.M. Mahomed, Int. J. Theor. Phys. \textbf{34,}
2267 (1995)

\bibitem{cla5} S.-F. Shen and Y.-Y. Jin, Acta. Math. Appl. Sinica (English
Series) \textbf{33,} 345 (2017)

\bibitem{clas6} S. Jamal, Gen. Relativ. Grav. \textbf{49,} 88 (2017)

\bibitem{leach1} P.G.L. Leach, J. Phys. A: Math. Gen. \textbf{13,} 1991
(1980)

\bibitem{leach2} S.\ Moyo and P.G.L. Leach,J. Math. Anal. Appl.\textbf{\ 252,%
} 840 (2000)

\bibitem{leach3} W. Sarlet, F.M. Mahomed and P.G.L. Leach, \ J. Phys. A.:
Math. Gen. \textbf{20,} 277 (1987)

\bibitem{leach4} F.M. Mahomed and P.G.L. Leach, Quaestiones Mathematicae 
\textbf{12,} 121 (1989)

\bibitem{erm} M. Tsamparlis and A. Paliathanasis, J. Phys. A: Math. Theor. 
\textbf{45,} 275202 (2012)

\bibitem{katz1} G.H.\ Katzin and J. Levine, J. Math. Phys. \textbf{9,} 8
(1968)

\bibitem{katz2} G.H.\ Katzin and J. Levine, J. Math. Phys. \textbf{22,} 1878
(1981)

\bibitem{katz3} G.H. Katzin, J. Math. Phys. \ \textbf{14,} 1213 (1973)

\bibitem{ami1} A.V. Aminova, Sbornik Math. \textbf{186}, 1711 (1995)

\bibitem{ami2} A.V. Aminova, Tensor N.S., \textbf{65}, 62 (2000)

\bibitem{p1} M. Tsamparlis and A. Paliathanasis, Nonlinear Dyn. \textbf{62,}
203 (2010)

\bibitem{p2} M. Tsamparlis and A. Paliathanasis, Gen. Relativ. Gravit. 
\textbf{42,} 2957 (2010)

\bibitem{p3} A. Paliathanasis and M. Tsamparlis, J. Geom. Phys. \textbf{62,}
2443 (2012)

\bibitem{p4} A. Paliathanasis and M. Tsamparlis, J. Geom. Phys. \textbf{107}%
, 45 (2016)

\bibitem{pp05} Y. Bozhkov and I.L. Freire, J. Diff. Eq. \textbf{249}, 872
(2010)

\bibitem{p6} A. Paliathanasis and M. Tsamparlis, Phys. Rev. D. \textbf{90},
043529(2014)

\bibitem{p7} A. Paliathanasis, M. Tsamparlis and M.T. Mustafa, Int. J. Geom.
Meth. Mod. Phys. \textbf{12}, 1550033 (2015)

\bibitem{p8} A. Paliathanasis, M. Tsamparlis and M.T. Mustafa, Comm. Non.
Sci. Num. Sim. \textbf{55}, 68 (2018)

\bibitem{chris} T.Christodoulakis, N. Dimakis and P.A. Terzis, J. Phys. A:
Math. Theor. \textbf{47}, 095202 (2014)

\bibitem{chris2} T. Christodoulakis, N. Dimakis, P.A. Terzis, Th. Grammenos,
E. Melas and A. Spanou, J. Geom. Phys.\textbf{\ 71}, 127 (2013)

\bibitem{chris3} T. Christodoulakis, N. Dimakis, P.A. Terzis and G. Doulis,
Phys. Rev.\ D \textbf{90}, 024052 (2014)

\bibitem{chris4} N. Dimakis, A. Karagiorgos, T. Pailas, P.A. Terzis and T.
Christodoulakis, Phys. Rev. D \textbf{95}, 086016 (2017)

\bibitem{anp1} A. Paliathanasis, M. Tsamparlis, S. Basilakos and J.D.
Barrow, Phys. Rev.\ D \textbf{93}, 043528 (2016)

\bibitem{nuz1} M.C. Nucci, Theor. Math. Phys. \textbf{168,} 994 (2011)

\bibitem{nuz2} M.C. Nucci, J. Phys.: Conf. Ser. \textbf{482}, 012032 (2014)

\bibitem{nuz4} M.C. Nucci and P.G.L. Leach, J. Non. Math. Phys.\textbf{\ 17}%
, 485 (2010)

\bibitem{pp1} B. Abraham-Shrauner, K.S. Govinder and D.J. Arrigo, J. Phys.
A: Math. Gen.\textbf{\ 39}, 5739 (2006)

\bibitem{pp2} M.L. Gandarias, J. Math. Anal. Appl. \textbf{348}, 752 (2008)

\bibitem{pp3} M. Tsamparlis and A. Paliathanasis, J. Geom. Phys. \textbf{73}%
, 209 (2013)
\end{thebibliography}
\end{document}